\documentclass[11pt,twoside,a4paper]{article}
\usepackage{amsmath,amssymb,amsthm}
\usepackage{graphicx}
\usepackage[english]{babel}
\usepackage{scrlayer-scrpage}
\usepackage[explicit]{titlesec}
\usepackage{libertine}
\usepackage{tikz}
\usetikzlibrary{calc}
\usepackage{tikzpagenodes}
\usepackage{hyperref}
\usepackage[margin=1.5in]{geometry}

\theoremstyle{plain}
\newtheorem{St}{Theorem}[section]
\newtheorem{Le}[St]{Lemma}
\newtheorem{Gev}[St]{Corollary}
\theoremstyle{definition}
\newtheorem{Ex}[St]{Example}
\newtheorem{Def}[St]{Definition}
\newtheorem{Opm}[St]{Remark}

\DeclareMathOperator\pg{\mathrm{PG}}
\DeclareMathOperator\ag{\mathrm{AG}}

\newcommand{\cL}{\mathcal{L}}
\newcommand{\gauss}[2]{\genfrac{[}{]}{0pt}{}{#1}{#2}}
\begin{document}
	\title{On two non-existence results for Cameron-Liebler $k$-sets in $\pg(n,q)$}

 \footnotetext{$*$ Department of Mathematics: Analysis, Logic and Discrete Mathematics, Ghent University, Krijgslaan 281, Building S8, 9000 Gent, Belgium (Email: leo.storme@ugent.be). \newline (http://cage.ugent.be/$\sim$ls).}%
    \footnotetext{$\dagger$ Department of Mathematics and Data Science, Vrije Universiteit Brussel,  Pleinlaan 2, 1050 Brussels, Belgium  (Email: Jan.De.Beule@vub.be, Jonathan.Mannaert@vub.be).}%

\author{Jan De Beule$^\dagger$, Jonathan Mannaert$^\dagger$ and Leo Storme$^*$}
\maketitle

\bibliographystyle{plain}
\begin{abstract}
This paper focuses on non-existence results for Cameron-Liebler $k$-sets.
A Cameron-Liebler $k$-set is a collection of $k$-spaces in $\pg(n,q)$ or $\ag(n,q)$ admitting a 
certain parameter $x$, which is dependent on the size of this collection. 
One of the main research questions remains the (non-)existence of 
Cameron-Liebler $k$-sets with parameter $x$. This paper improves two non-existence results. 
First we show that the parameter of a non-trivial Cameron-Liebler $k$-set in $\pg(n,q)$ should be 
larger than $q^{n-\frac{5k}{2}-1}$, which is an improvement of an earlier known lower bound. 
Secondly, we prove a modular equality on the parameter $x$ of Cameron-Liebler $k$-sets in $\pg(n,q)$ 
with $x<\frac{q^{n-k}-1}{q^{k+1}-1}$, $n\geq 2k+1$, $n-k+1\geq 7$ and $n-k$ even. In the affine
case we show a similar result for $n-k+1\geq 3$ and $n-k$ even. This is a generalization of earlier 
known modular equalities in the projective and affine case.

\end{abstract}
\section{Introduction}

Cameron and Liebler studied in \cite{Cameron-Liebler} irreducible collineation groups of $\pg(d,q)$ having 
equally many point orbits as line orbits. Such a group gives rise to a symmetrical tactical decomposition on $\pg(d,q)$, and
any line class in such a tactical decomposition is called a Cameron-Liebler line class. For $d=3$, a 
Cameron-Liebler line class is characterized by the property that it meets any spread of $\pg(3,q)$ 
in a fixed number $x$ of lines, where a spread of $\pg(3,q)$ is a partition of the point set in lines. 
So the parameter $x$ depends only on the size of such a class.

For $d = 3$, it is easy to see that the following line sets are examples of Cameron-Liebler line classes: 
(1) the empty set, with parameter $x=0$, (2) 
all lines through a fixed point $p$, with parameter $x=1$, (3) all lines in a fixed plane $\pi$, also with 
parameter $x=1$, and (4) a  union of (2) and (3) with $p\not\in\pi$, with parameter $x=2$. 
These examples (and their complements in the set of lines) are called {\em trivial examples}, 
and it was conjectured in \cite{Cameron-Liebler} that no other examples exist. This conjecture
has been disproven by Drudge in \cite{Drudge}, who gave an example in $\pg(3,3)$ with parameter $x=5$; 
an example that was generalized to an infinite family with parameter $x=\frac{q^2+1}{2}$ in $\pg(3,q)$, $q$ odd, 
in \cite{BruenAndDrude}. New non-trivial examples have been discovered by Rodgers in \cite{RodgersPhD,Rodgers2013}, 
some of them have been generalized to infinite families, see \cite{DeBeule2016,Feng20xx,Feng2015}. 
Generally spoken, non-trivial examples are rare. Furthermore, non-existence results of Cameron-Liebler 
line classes for particular values of the parameter $x$  have been found, see e.g. 
\cite{MetschAndGavrilyuk}. Non-trivial examples remain rare for $d=3$, and there has also been 
much effort to show non-existence results. For $d=3$, probably one of the most powerful non-existence
results is found in \cite{MetschAndGavrilyuk}, excluding roughly half of the possible parameters. 
Cameron-Liebler $k$-sets occur in applications in coding theory. A Cameron-Liebler line 
class in $\pg(3,q)$ can be described as a complete regular code in the Grassmann graph that
is not a $t$-design. These codes were studied earlier in \cite{Delsarte1973}. 
Recent work in this area is e.g. \cite{Mogilnykh2022,SantZullo}.

Cameron-Liebler sets of $k$-spaces in $\pg(n,q)$ or $\ag(n,q)$ have also been studied quite intensively. We refer
to the following incomplete list \cite{Jozefien, Me1, Me2, DrudgeThesis}. So far, most non-existence results 
for $n \geq 3$ are either formulated as a lower bound on the parameter $x$. Recently, in \cite{Me4}, a modular equality
on the parameter $x$ of Cameron-Liebler $k$-sets was shown comparable to the one found in \cite{MetschAndGavrilyuk}, albeit slightly
weaker. Very recently, an asymptotic non-existence result for non-trivial Cameron-Liebler $k$-sets in general dimension $n > 3$ 
was shown in \cite{Ihringer2023}, albeit for $n$ {\em much larger than 3}. The details are
too technical to capture in this introduction. Although this very recent result might indicate that 
non-trivial Cameron-Liebler $k$-sets in dimension $n > 3$ are unexpected to exist, there is 
still room to look for examples. Furthermore, the result from \cite{Ihringer2023} still 
leaves room to improve on non-existence conditions on the parameter $x$. 

In this paper, we will focus on Cameron-Liebler sets of $k$-spaces (or CL $k$-sets for short) in $\pg(n,q)$ 
and $\ag(n,q)$, where, in general, we choose $n\geq 2k+1\geq 3$. 

To define Cameron-Liebler $k$-sets in $\pg(n,q)$, respectively $\ag(n,q)$, we use the point-($k$-space) incidence matrix,
denoted by $P_n$, respectively $A_n$. This matrix is a $\{0,1\}$-valued matrix, rows indexed 
by the points, and columns indexed by the $k$-spaces, with $0$ if the point is not incident with the $k$-space 
and $1$ otherwise. 

The \emph{Gaussian binomial coefficient} is defined as 
$$\gauss{n}{k}_q:=\frac{(q^n-1)\cdots (q^{n-k+1}-1)}{(q^k-1)\cdots (q-1)}$$
and equals the number of $(k-1)$-spaces inside $\pg(n-1,q)$.

\subsection{Cameron-Liebler $k$-sets in $\pg(n,q)$}

\begin{Def}\label{th:equivdef}
A Cameron-Liebler set of $k$-spaces (or CL $k$-set) is a collection $\mathcal{L}$  of $k$-spaces in $\pg(n,q)$, for which its characteristic vector
$$\chi_\mathcal{L}\in \mathrm{Im}(P_n^T).$$
The set $\mathcal{L}$ has parameter $x$ if and only if 
$$|\mathcal{L}|=x\gauss{n}{k}_q.$$
\end{Def}
From the definition, it follows that $0\leq x\leq \frac{q^{n+1}-1}{q^{k+1}-1}$. We list a number of (well known) equivalent definitions for CL $k$-sets 
 in the following theorem.
\begin{St}\label{th:EquivalenceProj}\label{EquivalenceProj}\cite[Theorem 2.2]{Jozefien}
	Let \(\mathcal{L}\) be a non-empty set of $k$-spaces in PG$(n,q)$, $n \geq 2k+1$, and $x$ so that 
	$|\mathcal{L}|= x\gauss{n}{k}_q$. Then the following properties are equivalent.
	\begin{enumerate}
		\item \(\mathcal{L}\) is a Cameron-Liebler $k$-set in $\pg(n,q)$.
		\item For every $k$-space $K$, the number of elements of $\mathcal{L}$ disjoint from $K$ is equal to $(x-\chi_\mathcal{L}(K))\gauss{n-k-1}{k}_q q^{k^2+k}$.
\item For an $i \in \{1,...,k+1\}$ and a given $k$-space $K$, the number of elements of $\mathcal{L}$, meeting $K$ in a $(k-i)$-space is given by:
		$$\left\{\begin{array}{ll}
			  \left( (x-1) \frac{q^{k+1}-1}{q^{k-i+1}-1}+ q^i \frac{q^{n-k}-1}{q^i-1}\right) q^{i(i-1)} \begin{bmatrix}
			n-k-1 \\
			i-1
			\end{bmatrix}_q \begin{bmatrix}
			k \\
			i
			\end{bmatrix}_q & \text{ if } K \in \mathcal{L} \\
			xq^{i(i-1)} \begin{bmatrix}
			n-k-1 \\
			i-1
			\end{bmatrix}_q \begin{bmatrix}
			k+1 \\
			i
			\end{bmatrix}_q & \text{ if } K \not\in \mathcal{L}.
		\end{array} \right.$$
	
		\item If $(k+1) \mid (n+1)$, i.e. if and only if PG$(n,q)$ has $k$-spreads, then $|\mathcal{L}\cap \mathcal{S}|=x$ for any $k$-spread $\mathcal{S}$. 
	\end{enumerate}
\end{St}

From Theorem \ref{th:EquivalenceProj} (4), it follows that $x$ is always an integer if $(k+1)\mid (n+1)$. 

For $n > 3$, the trivial examples are the natural generalizations of the trivial examples of a CL line class: (1) the empty set ($x=0$), 
(2) all $k$-spaces through a point ($x=1$), (3) all $k$-spaces contained in a hyperplane ($x=\frac{q^{n-k}-1}{q^{k+1}-1}$), 
and (4) the union of (2) and (3) if the point is not lying in the hyperplane ($x=1+\frac{q^{n-k}-1}{q^{k+1}-1}$). 
All complements of these examples are also CL $k$-sets, and all examples from this list are called {\em trivial}. 
A general result to be used will be the following lemma. 

\begin{Le}\cite[Lemma 2.12]{Jozefien}\label{DrudgeArgum}
Let $\mathcal{L}$ be a Cameron-Liebler $k$-set in $\pg(n,q)$, then we find the 
following equality for every point $p$ and every $i$-dimensional subspace $\tau$, 
with $p\in \tau$ and $i\geq k+1$,
\[
\left| [p]_k\cap \mathcal{L} \right|+ \frac{\gauss{n-1}{k}_q (q^k-1)}{\gauss{i-1}{k}_q (q^i-1)} \left| [\tau]_k\cap \mathcal{L} \right| 
= \frac{\gauss{n-1}{k}_q }{\gauss{i-1}{k}_q } \left| [p,\tau]_k\cap \mathcal{L} \right| + \frac{q^k-1}{q^n-1} \left| \mathcal{L} \right|\,.
\]
Here $[p]_k, [\tau]_k$ and $[p,\tau]_k$ denote all $k$-spaces of $\pg(n,q)$ containing the point $p$, contained in $\tau$ or both respectively.
\end{Le}

\subsection{Cameron-Liebler $k$-sets in $\ag(n,q)$}
Similarly, we can define a CL $k$-set in $\ag(n,q)$.  
\begin{Def}
A Cameron-Liebler set of $k$-spaces (or CL $k$-set) is a collection $\mathcal{L}$ of $k$-spaces in $\ag(n,q)$, for which its characteristic vector
$$\chi_\mathcal{L}\in \text{Im}(A_n^T).$$
The CL $k$-set  $\mathcal{L}$ has parameter $x$ if and only if 
$$|\mathcal{L}|=x\gauss{n}{k}_q.$$
\end{Def}

We list a number of (well known) equivalent definitions for CL $k$-sets in affine spaces in the following theorem.
\begin{St}\cite[Theorem 3.5]{Me2}\label{th:EquivalenceAffine}\label{2to8}
	Consider the affine space $\ag(n,q)$, for $n\geq 2k+1$, and let \(\mathcal{L}\) be a set of $k$-spaces such that $|\mathcal{L}|=x$\mbox{$\scriptsize{\begin{bmatrix}
	n \\
	k
	\end{bmatrix}_q}$} for a positive integer $x$. Then the following properties are equivalent.
	\begin{enumerate}
		\item \(\mathcal{L}\) is a Cameron-Liebler $k$-set in $\ag(n,q)$.
		\item For every $k$-spread $\mathcal{S}$, it holds that $|\mathcal{L}\cap \mathcal{S}|=x.$
		\item For every pair of conjugated switching $k$-sets \(\mathcal{R}\) and \(\mathcal{R}'\), \(|\mathcal{L} \cap \mathcal{R} |= |\mathcal{L} \cap \mathcal{R}'|\).
\end{enumerate}
If $k=1$ and we thus consider Cameron-Liebler line classes, then the following property is equivalent to the previous ones.
\begin{enumerate}
          \item[4.] For every line $\ell$, the number of elements of \(\mathcal{L}\) affinely disjoint to $\ell$ is equal to
	\begin{equation} \label{EqLinesEq} \left( q^2\begin{bmatrix}
	n-2\\
	1
	\end{bmatrix}_q +1 \right) (x- \chi_{\mathcal{L}}(\ell))\end{equation}
	and through every point at infinity there are exactly $x$ lines of \(\mathcal{L}\).

	\end{enumerate}
\end{St}

For CL $k$-sets in $\ag(n,q)$, the parameter $x$ is always an integer and $0\leq x\leq q^{n-k}$.\\
Once again, examples of CL $k$-sets in $\ag(n,q)$ are rare. For $n > 3$, no other examples than the following are known so far.
These are: (1) the empty set ($x=0$), (2) all $k$-spaces through a point ($x=1$), and their complements. These examples
are called {\em trivial}. Furthermore, most of the known examples of CL line classes in $\pg(3,q)$ turn out to be a CL line class
in $\ag(3,q)$. For more details on CL $k$-sets in $\ag(n,q)$, we refer to \cite{Me2}.

\subsection{Basic properties}

The following lemma summarizes elementary properties of CL $k$-sets that follow from Theorem \ref{th:EquivalenceProj} and Theorem \ref{th:EquivalenceAffine}. 


\begin{Le}\cite[Lemma 3.1]{Jozefien},\cite[Lemma 26]{Me2}
Suppose that $\mathcal{L}$ and $\mathcal{L}'$ are two Cameron-Liebler $k$-sets both in $\pg(n,q)$, respectively $\ag(n,q)$, with parameter $x$ and $x'$. Then the following properties
hold.
\begin{itemize}
\item The complement of $\mathcal{L}$ is a Cameron-Liebler $k$-set with parameter $\frac{q^{n+1}-1}{q^{k+1}-1}-x$, respectively $q^{n-k}-x$.
\item If $\mathcal{L}$ and $\mathcal{L}'$ are disjoint as sets of $k$-spaces, then $\mathcal{L}\cup \mathcal{L}'$ is a Cameron-Liebler $k$-set with parameter $x+x'$.
\item If $\mathcal{L}\subseteq \mathcal{L}'$, then $\mathcal{L}'\setminus \mathcal{L}$ is a Cameron-Liebler $k$-set with parameter $x'-x$.
\end{itemize}
\end{Le}
The first statement of this lemma implies that any classification result for the first half of the parameters yields a full classification. 
The following property can be shown relatively easily, and explains how CL line classes in $\pg(3,q)$ can be CL line classes in $\ag(3,q)$
as well.

\begin{St}\cite[Theorem 1]{Me2}\label{th:ProjToAff}\label{H6ProjToAff}
Let $\cL$ be a Cameron-Liebler $k$-set with parameter $x$ in PG($n,q$) which does not contain $k$-spaces in some hyperplane $H$. 
Then $\cL$ is a Cameron-Liebler $k$-set with parameter $x$ of $\ag(n, q) = \pg(n, q) \setminus H$.
\end{St}
Conversely, a similar result holds.
\begin{St}\cite[Theorem 2]{Me2}\label{th:affToProj}\label{CLkSpaceBasic}
    If $\cL$ is a Cameron-Liebler $k$-set of $\ag(n, q)$ with parameter $x$, then $\mathcal{L}$ is a Cameron-Liebler $k$-set of $\pg(n, q)$
    with parameter $x$ in the projective closure $\pg(n, q)$ of $\ag(n, q)$.
\end{St}

\section{Preliminaries}

In this section, we summarize non-existence results for CL $k$-sets with certain parameters. Such results are either a lower bound on the parameter $x$
or a modular equality on $x$. 

\subsection{Non-existence results in $\pg(n,q)$}

Classification or non-existence results for CL $k$-sets with small parameters are summarized in the following theorems.
\begin{St}\cite[Theorem 4.3]{Jozefien}\label{th:Non-existence01}\label{Non-existence01}
There do not exist CL $k$-sets in $\pg(n,q)$ with parameter $x \in ]0,1[$ and if $n\geq 3k+2$, then there are no CL $k$-sets with parameter $x\in]1,2[$.
\end{St}

\begin{St}\cite[Theorem 4.1]{Jozefien}\label{th:Uniquenessx=1}\label{Uniquenessx=1}
Let $\mathcal{L}$ be a CL $k$-set  with parameter $x = 1$ in $\pg(n,q)$, $n \geq 2k+1$. Then $\mathcal{L}$ consists out of all the $k$-spaces through a fixed point, or $n = 2k+1$ and $\mathcal{L}$ is the set of all the $k$-spaces in a hyperplane of $\pg(2k + 1,q)$.
\end{St}

The following theorem provides a lower bound on the parameter $x$.

\begin{St}\cite[Theorem 1.1]{Me3}\label{th:MainOld}\label{MainOld}
Suppose that $n\geq 3k+3$ and $k\geq 1$. Let $\mathcal{L}$ be a CL $k$-set with parameter $x$ in $\pg(n,q)$  such that $\mathcal{L}$ is not a point-pencil, nor the empty set. Then 
$$x\geq\frac{q^{n-k}-1}{q^{2k+2}-1}+1.$$
\end{St}
\begin{Opm}
Note that the lower bound of Theorem \ref{th:MainOld} is roughly $q^{n-3k-2}$.
\end{Opm}
This was very recently improved, resulting in the following theorem.
\begin{St}\label{th:boundFerd}\cite[Theorem 5.1]{Ihringer2023}
    Suppose that $n\geq 3k$ and $k\geq 1$. Let $\mathcal{L}$ be a CL $k$-set with parameter $x$ in $\pg(n,q)$  such that $\mathcal{L}$ is not a point-pencil, nor the empty set. Then 
$$x\geq \frac{(q^n-1)(q-1)^2}{q(q^k-1)^2(q^{k+1}-1)}>\frac{1}{4}q^{n-3k}.$$
\end{St}

The following results are modular equalities on the parameter $x$. In combination with the bounds above they show to be very useful in excluding parameters.

\begin{St}\cite[Theorem 1.1]{MetschAndGavrilyuk}\label{th:GMModEq}\label{MetschRes}
Suppose that \(\mathcal{L}\) is a CL line class with parameter \(x\) of $\pg(3,q)$.  
Then for every plane and every point of $\pg(3,q)$, 
\begin{equation}\label{Metsch}
\binom{x}{2} +m(m-x)\equiv 0 \mod (q+1),
\end{equation}
where \(m\) is the number of lines of \(\mathcal{L}\) in the plane, respectively through the point.
\end{St} 

A slightly weaker result was shown for $n \geq 7$ odd in \cite{Me4}. 

\begin{St}\cite[Theorem 1.3]{Me4}\label{th:projectiveLines}
Suppose that $\mathcal{L}$ is a CL line class with parameter $x$ in $\pg(n,q)$, $n\geq 7$ odd. Then for
any point $p$,
\[
x(x-1)+2\overline{m}(\overline{m}-x)\equiv 0 \mod (q+1)\,,
\]
where $\overline{m}$ is the number of lines of $\mathcal{L}$ through $p$.
\end{St}

Finally, the following result classifies all CL $k$-sets asymptotically in large dimensions.  
\begin{St}\label{th:ihringer}\cite[Theorem 1.1]{Ihringer2023}
For each $k \geq 1$ and $q$, there exists a natural number $c_0(k,q)$ such that all CL $k$-sets in $\pg(n,q)$ are trivial 
if $\max\{k,n-k\}\geq c_0(k,q).$
\end{St}
As indicated clearly in \cite{Ihringer2023},  this rules out CL $k$-sets only in very large dimensions. 


\subsection{Non-existence result in $\ag(n,q)$}

Using Theorem \ref{th:affToProj}, it is easy to see that all results for $\pg(n,q)$ are still valid. However, since examples in $\ag(n,q)$ are examples in the projective closure of a very specific form, these results can often be improved. This improvement can be seen in the following results.

\begin{St}\cite[Theorem 6.5]{Me3}\label{th:NonExAG}\label{NonExAG}
Suppose that $n\geq 2k+2$ and $k\geq 1$. Let $\mathcal{L}$ be a CL $k$-set of parameter $x$ in $\ag(n,q)$ such that $\mathcal{L}$ is not a point-pencil, nor the empty set. Then 
$$x\geq 2\left(\frac{q^{n-k}-1}{q^{k+1}-1}\right)+1.$$
\end{St}
\begin{Opm}
Here the lower bound of Theorem \ref{th:NonExAG} is roughly $q^{n-2k-1}$.
\end{Opm}

\begin{St}\cite[Theorem 1.5]{Me4}\label{th:affine-init}
Suppose that $\mathcal{L}$ is a CL line class in $\ag(n,q)$, $n \geq 3$ odd, with parameter $x$, then
$$x(x-1) \equiv0 \mod 2(q+1).$$
\end{St}

\subsection{Structure of the paper}
The paper has  two major sections.  In Section 3, we will improve Theorem \ref{th:boundFerd}, 
for $n\geq  3k+4$, by some inductive arguments.  However this will be done in several steps. First in Section 3.1, 
we only obtain a small improvement for $k>3$. Following in Section 3.2, we will improve these results even further 
and for general $k$. Both sections are required since Section 3.2 requires Section 3.1. Very atypical for this 
result is that the improvement is obtained by looking at some underlying affine spaces. Finally, in Section 4, 
we will discuss the generalization of the modular equality proven in Theorems \ref{th:projectiveLines} 
and \ref{th:affine-init} to sets of $k$-spaces.

\section{Improving the lower bound on the parameter $x$}

The goal of this section is to improve the lower bound on the parameter $x$ of non-trivial CL $k$-sets in $\pg(n,q)$. In particular, we aim to improve Theorem \ref{th:ihringer}. In order to do this, we require some introductory results. The techniques we use in this section are similar to those in \cite{Me3}. Nevertheless, it can be interesting to compare both. For this we refer to \cite[Theorem 1.1]{Me3}.

\begin{St}[Folklore]\label{th:Subsp}\label{Subsp}
Consider a CL $k$-set $\mathcal{L}$ in $\pg(n,q)$, with $n \geq 2k+1$, 
and let $\pi$ be a subspace of dimension $i\geq k+1$. Then $\mathcal{L}_\pi:=\{K\in \mathcal{L} \mid K\subseteq \pi\}$  is also a CL $k$-set in $\pi$.
\end{St}

The relation between the parameter $x$ of $\mathcal{L}$ and the parameter of the induced CL $k$-sets $\mathcal{L}_\pi$ is 
given by the following theorem. 


\begin{Le}\cite[Lemma 4.1]{Me3}\label{Le:Aid}\label{Aid}
Suppose that $\mathcal{L}$ is a CL $k$-set in $\pg(n,q)$ with parameter $x$. 
Then for every  $t$, such that  $2k+1 \leq t\leq n-1$ and $n\geq 2k+2$, we can fix an arbitrary $k$-space $K$ in $\pg(n,q)$. Consider all $t$-dimensional subspaces $\pi_i$ through $K$, each admitting a CL $k$-set $\mathcal{L}_{\pi_i}$ with parameter $x_{\pi_i}$. Then we find that
\begin{equation}
x=\frac{\left(\sum_{K \subseteq \pi_i}x_{\pi_i}-  \gauss{n-k}{t-k}_q\chi_\mathcal{L}(K) \right)} {\gauss{n-k-1}{t-k-1}_q }+\chi_\mathcal{L}(K)\,,
\end{equation}
where  $\chi_{\mathcal{L}}(K)$ is the value of the characteristic vector of $\mathcal{L}$ at position $K$
and where the sum runs over all $t$-spaces $\pi_i$ through $K$.
\end{Le}

The following theorem describes how a CL $k$-set can be shown to be trivial if its induced CL $k$-set with relation to a particular 
subspace is trivial.

\begin{St}\cite[Theorem 3.7]{Me3} \label{th:DrudgeGeneralPG}
Let $n\geq 2k+1$. Suppose that $\mathcal{L}$ is a CL $k$-set in $\pg(n,q)$. 
Suppose that there exists an $i$-space $\pi$, with $i\geq k+1$, such that $\mathcal{L}$ restricted to $\pi$ consists of all $k$-spaces through a point $p$.
Then $\mathcal{L}$ is the  set of $k$-spaces 
through this same point $p$. 
\end{St}

Now we can prove the following theorem.

\begin{St}\label{Le:Aid2}\label{Aid2}
Suppose that for $t\geq 3k+2$, there exists a lower bound $B(t,q)$ such that for every CL $k$-set with parameter $x$ in $\pg(t,q)$ either (1) $x\in\{0,1\}$, or (2) $x\geq B(t,q)$. Let $n\geq t+1$, then the parameter $x$ of a Cameron-Liebler $k$-set in $\pg(n,q)$ satisfies (1) $x\in\{0,1\}$, or (2) $x\geq \left(B(t,q)-1\right)\frac{q^{n-k}-1}{q^{t-k}-1}+1$.
\end{St}
\begin{proof}
Consider an arbitrary fixed $k$-space $K\in \mathcal{L}$. Here we assume that $\mathcal{L}\not= \emptyset$. 
Consider all $t$-spaces $\pi_i$ containing $K$, each inducing a different CL $k$-set in every $t$-space. Using Lemma~\ref{Le:Aid}, we find that
\begin{equation*}
x=\frac{\left(\sum_{K \subseteq \pi_i}x_{\pi_i}-  \gauss{n-k}{t-k}_q\chi_\mathcal{L}(K) \right)} {\gauss{n-k-1}{t-k-1}_q }+\chi_\mathcal{L}(K)\,.
\end{equation*}
Since $K\in \mathcal{L}$ in every induced CL $k$-set, we find that $x_{\pi_i} \not= 0$. Using Theorem \ref{th:Non-existence01}, we can improve this to $x_{\pi_i} \geq 1$. This can further be improved by assuming that $\mathcal{L}$ is not the collection of all $k$-spaces containing a fixed point, where from Theorem \ref{th:Uniquenessx=1} and Theorem \ref{th:DrudgeGeneralPG}, it would follow that $x_{\pi_i} >1$. Finally, we can use the assumption to obtain that $x_{\pi_i} \geq B(t,q)$. After substituting this in the previous equation, we obtain the assertion.
\end{proof}

\begin{Opm}
Note that Theorem \ref{th:MainOld} is in fact a special case of Theorem \ref{Le:Aid2}, where, due to Theorem \ref{th:Non-existence01} and Theorem \ref{th:Uniquenessx=1}, it was obtained that $B(3t+2,q)>2$.
\end{Opm}

\subsection{A first improvement of the lower bound}
First we take a look at the following non-existence result.

\begin{St}\cite[Theorem 4.9]{Jozefiencorr}\label{th:NonExProjOld}\label{JozefienClassific}
There are no CL $k$-sets in $\pg(n,q)$, with $n\geq 3k+2$ and $q\geq 3$, with parameter $x$ if
	$$2\leq x\leq C(n,k,q),$$
with $C(n,k,q)=\frac{1}{\sqrt[8]{2}}q^{\frac{n}{2}-\frac{k^2}{4}-\frac{3k}{4}-\frac{3}{2}}(q-1)^{\frac{k^2}{4}-\frac{k}{4}+\frac{1}{2}}\sqrt{q^2+q+1}$.
\end{St}
\begin{Opm}
Note that in \cite[Theorem 4.1]{BDF}, it was shown that all CL sets of $k$-spaces in $\pg(n,q)$, 
with $n\geq 2k+1\geq 5$, $k\geq 2$ and $q\in\{2,3,4,5\}$, are trivial. Hence,
Theorem \ref{th:NonExProjOld} also holds for $q=2$, as was pointed out by Ferdinand Ihringer \cite{private}.
\end{Opm}

Using this theorem, we obtain the following result.
\begin{St}\label{th:MainNew}\label{MainNew}
Suppose that $n\geq 3k+3$ and $k\geq 1$. Let $\mathcal{L}$ be a CL $k$-set with parameter $x$ in $\pg(n,q)$  such that $\mathcal{L}$ is neither a point-pencil, nor the empty set. Then 
$$x\geq D(k,q)\left(\frac{q^{n-k}-1}{q^{2k+2}-1}\right)+1,$$
for $D(k,q) := C(3k+2,k,q)=\frac{1}{\sqrt[8]{2}}q^{\frac{3k}{4}-\frac{k^2}{4}-\frac{1}{2}}(q-1)^{\frac{k^2}{4}-\frac{k}{4}+\frac{1}{2}}\sqrt{q^2+q+1}-1$.
\end{St}
\begin{proof}
Use Theorem \ref{Le:Aid2} in combination with Theorem \ref{th:NonExProjOld}. Choose $t=3k+2$ for the optimal result.
\end{proof}

\begin{Opm}
We can see that in Theorem \ref{th:MainNew}, the bound of Theorem \ref{th:MainOld} is 
multiplied by a factor $D(k,q)$, which itself is of size roughly 
$q^{\frac{k}{2}+1}$. Hence, the lower bound will now be of size 
roughly $q^{n-\frac{5k}{2}-1}$.
\end{Opm}

\begin{Opm}[recursive argument]
One might have expected that Theorem \ref{th:MainOld} would give an improvement compared with using Theorem 
\ref{th:NonExProjOld} in the proof of Theorem~\ref{th:MainNew}, 
since the first is a stronger non-existence result. 
However, this is not the case, as illustrated by the following argument. 
Consider $t\geq 3k+3$ and $n\geq t+1$. Hence, we obtain that
\begin{equation*}
\begin{split}
B(n,q)&=(B(t,q)-1)\frac{q^{n-k}-1}{q^{t-k}-1}+1\\
&= \left(\frac{q^{t-k}-1}{q^{2k+2}-1} +1-1\right)\frac{q^{n-k}-1}{q^{t-k}-1}+1=\frac{q^{n-k}-1}{q^{2k+2}-1}+1,
\end{split}
\end{equation*}
which does not yield a better bound. Similarly, a worse result is obtained by filling in Theorem \ref{th:boundFerd} for $t=3k+2$. From this bound, we obtain roughly $x\geq q^{n-3k}$.
\end{Opm}

\subsection{A second improvement of the lower bound}

We now proceed in two steps. The bounds on $x$ in Theorem~\ref{th:MainAff} and \ref{th:MainNewNew} are equal. 
But Theorem~\ref{th:MainAff} needs the extra condition that every hyperplane contains at least one element
of $\cL$. We have to work in two steps to remove the condition in Theorem~\ref{th:MainNewNew}.

\begin{St}\label{th:MainAff}
Suppose that $n\geq 3k+4$ and $k\geq 1$. Consider a CL $k$-set $\cL$ with parameter $x$ in $\pg(n,q)$, which is neither empty, nor a point-pencil. 
If every hyperplane contains at least one element of $\cL$, then
$$x\geq \left(D(k,q)\frac{q^{n-k-1}-1}{q^{2k+2}-1}+1\right)\frac{q^{n+1}-1}{q^{n}-1}\,,$$
with $D(k,q) = \frac{1}{\sqrt[8]{2}}q^{\frac{3k}{4}-\frac{k^2}{4}-\frac{1}{2}}(q-1)^{\frac{k^2}{4}-\frac{k}{4}+\frac{1}{2}}\sqrt{q^2+q+1}-1$.
\end{St}
\begin{proof}
Double count the pairs $(K, \pi)$, with $\pi$ a $t$-dimensional subspace and $K\in \cL$ such that $K\subseteq \pi$. 
\begin{enumerate}
\item If we fix a $k$-space $K\in \mathcal{L}$, then we know that there are $\gauss{n-k}{t-k}_q$ subspaces $\pi$ containing it. Consequently, the number of pairs equals
$$|\cL|\gauss{n-k}{t-k}_q.$$
\item Fixing an arbitrary $t$-space $\pi$, we know, by Theorem \ref{th:Subsp}, that $\cL_\pi:=\{K\in \mathcal{L} \mid K\subseteq \pi\}$ is a CL $k$-set in $\pi$ with a certain parameter $x_\pi$. Hence, we obtain that the number of pairs equals
$$\sum_{\pi\subseteq \pg(n,q)} |\cL_\pi|.$$
\end{enumerate}
Using that  $|\cL|=x\gauss{n}{k}_q$ and $|\cL_\pi|=x_\pi\gauss{t}{k}_q$, we find that
$$x\gauss{n}{k}_q\gauss{n-k}{t-k}_q=\gauss{t}{k}_q\sum_{\pi\subseteq \pg(n,q)} x_\pi.$$
Combining this with 
$$\frac{\gauss{t}{k}_q}{\gauss{n}{k}_q\gauss{n-k}{t-k}_q}=\frac{1}{\gauss{n}{t}_q},$$
results in
\begin{equation}\label{eq2}
x=\frac{\sum_{\pi\subseteqq \pg(n,q)} x_\pi}{\gauss{n}{t}_q}.
\end{equation}
Using $t=n-1$, in combination with the assumption that $\mathcal{L}$ as a set of $k$-spaces is not skew to any hyperplane, we obtain that for every hyperplane $\pi$, it holds that $x_{\pi}\not=0$. Similarly like in the proof of Theorem \ref{th:MainNew}, we may assume that $x_\pi>2$. Using the result of this theorem, we find that 
$$x_{\pi}\geq D(k,q)\frac{q^{n-k-1}-1}{q^{2k+2}-1}+1.$$ 

The assertion now follows by substituting this in Equation (\ref{eq2}).
\end{proof}

Often non-existence results for CL $k$-sets in $\ag(n,q)$ are based on results in $\pg(n,q)$. In the next theorem,
we will use results of CL $k$-sets in $\ag(n,q)$ to improve Theorem \ref{th:MainNew} in $\pg(n,q)$.

\begin{St}\label{th:MainNewNew}
Suppose that $n\geq 3k+4$ and $k\geq 1$. Let $\mathcal{L}$ be a CL $k$-set with parameter $x$ in $\pg(n,q)$  
such that $\mathcal{L}$ is neither a point-pencil, nor the empty set. Then 
$$x\geq  \left(D(k,q)\frac{q^{n-k-1}-1}{q^{2k+2}-1}+1\right)\frac{q^{n+1}-1}{q^{n}-1},$$
for $D(k,q)=\frac{1}{\sqrt[8]{2}}q^{\frac{3k}{4}-\frac{k^2}{4}-\frac{1}{2}}(q-1)^{\frac{k^2}{4}-\frac{k}{4}+\frac{1}{2}}\sqrt{q^2+q+1}-1$.
\end{St}
\begin{proof}
Using Theorem \ref{th:MainNew}, we already know that $x\geq  D(k,q)\left(\frac{q^{n-k}-1}{q^{2k+2}-1}\right)+1$. If $\mathcal{L}$ is a CL $k$-set in $\pg(n,q)$ with a parameter $x$ for which 
$$ D(k,q)\left(\frac{q^{n-k}-1}{q^{2k+2}-1}\right)+1\leq x < \left(D(k,q)\frac{q^{n-k-1}-1}{q^{2k+2}-1}+1\right)\frac{q^{n+1}-1}{q^{n}-1},$$
then, by Theorem \ref{th:MainAff}, we obtain that there exists at least one hyperplane $\pi$ 
not containing any element of $\mathcal{L}$.  In this case, we can use 
Theorem \ref{th:ProjToAff} to conclude that $\mathcal{L}$ defines a 
CL $k$-set with the same parameter $x$ in $\ag(n,q) = \pg(n,q) \setminus \pi$. Here 
the affine space is induced by the skew hyperplane. 
Using Theorem \ref{th:NonExAG} and the fact that
$$\left(D(k,q)\frac{q^{n-k-1}-1}{q^{2k+2}-1}+1\right)\frac{q^{n+1}-1}{q^{n}-1}<2\left(\frac{q^{n-k}-1}{q^{k+1}-1}\right)+1,$$
we obtain a contradiction.
\end{proof}

\subsection{Conclusion}
Theorem \ref{th:boundFerd} gives roughly $q^{n-3k}$ as a lower bound on $x$, 
while Theorem \ref{th:MainNew} yields roughly $q^{n-\frac{5k}{2}-1}$, 
which is an improvement whenever $k>2$. However, Theorem \ref{th:MainNewNew} improves 
this result only slightly by adding an additional 
term $q$, thus we obtain that $x$ is larger than roughly $q^{n-\frac{5k}{2}-1}+q$. 
Since every CL $k$-set in $\ag(n,q)$ is also a CL $k$-set in $\pg(n,q)$ with the 
same parameter $x$, we can compare with Theorem \ref{th:NonExAG}, but this gives
no further improvement. 

\begin{Ex}
Suppose that $q=7$, $k=1$ and $n=8$. Then we know that the parameter 
$x$ of a CL $k$-set in $\pg(n,q)$ lies in the interval 
$[0, \frac{q^{n+1}-1}{q^{k+1}-1}]=[0; 840700.125]$. In particular, 
we only need to classify those CL $k$-sets with parameter 
$x\leq 420350$, which is the first half of all integers of the range. 
Theorem \ref{th:MainNew} gives that $x\geq 5476.998...$, while 
Theorem \ref{th:MainNewNew} yields  that all parameters $x \geq 5482.9585...$.
Recall that the parameter $x$ is guaranteed to be an integer only if $(k+1) \mid (n+1)$, therefore
we did not round up the numerical values in this example.
\end{Ex}


\section{Modular equality}
\subsection{The projective case}
The goal of this section is to generalize the modular equality  in Theorem \ref{th:projectiveLines} given in \cite{Me4} to CL sets of $k$-spaces in $\pg(n,q)$ and to $\ag(n,q)$.
\begin{St}[Folklore]\label{DualSubspace}
Suppose that $\mathcal{L}$ is a CL $k$-set in $\pg(n,q)$, with $n\geq 2k+1$. 
Suppose that $\pi$ is an $i$-dimensional space for $i\leq k-1$, 
and let $\beta$ be any $(n-i-1)$-space skew to $\pi$. Then the set 
$\mathcal{L}^\pi:=\{K\cap \beta \mid K \in \mathcal{L}, \pi\subseteq K\}$ 
is a CL $(k-i-1)$-set in $\beta=\pg(n-i-1,q)$.
\end{St}
\begin{proof}
We will use induction on the dimension $i$. First assume that $i=0$. 
Suppose that $p$ is a fixed point and that $\beta$ is an $(n-1)$-space 
in $\pg(n,q)$ not containing $p$. We will show that if we project
all $k$-spaces of $\cL$ containing $p$ from $p$ on $\beta$, we obtain 
a CL $(k-1)$-set in $\beta$. Consider the point-($k$-space) incidence matrix $P_n$
of $\pg(n,q)$ and index the rows in such a way that the first row corresponds to the point
$p$, and the first columns are indexed by the $k$-spaces of $\pg(n,q)$ through $p$. 
Define $\mathcal{P}_1:=\{p\}$, $\mathcal{P}_{2,1}:=\{p \mid p\in \beta\}$ and 
$\mathcal{P}_{2,2}$ all other points. Similarly 
$\mathcal{B}_1$ is the set of all $k$-spaces through $p$, and 
$\mathcal{B}_2$ the set of all other $k$-spaces. The following
conclusions are clear.
%
\begin{itemize}
    \item The point-($k$-space) incidence matrix of $\mathcal{P}_1$ and $\mathcal{B}_1$ equals the transpose of the all-one vector, denoted by $\mathbf{1}^T$.
    \item Similarly, the point-($k$-space) incidence matrix of $\mathcal{P}_1$ and $\mathcal{B}_2$ equals the transpose of the all-zero vector, denoted by $\mathbf{0}^T$.
    \item Finally, to describe the point-($k$-space) incidence matrix of $\mathcal{P}_{2,1}\cup \mathcal{P}_{2,2}$ and $\mathcal{B}_1$ is
    is slightly harder. For every point $p_1\in \mathcal{P}_{2,1}$ it follows that a $k$-space $K$ through $p$ contains $p_1$ if and only if the projection of $K$ from $p$ onto $\beta$ contains $p_1$. This projection, in its turn is a $(k-1)$-space. Hence the point $p_1$ indexes
    a row that is also found in the point-$(k-1)$-space incidence matrix of $\pg(n-1,q)$. For points $p_2\in \mathcal{P}_{2,2}$, we consider the projection of $p_2$ from $p$ onto $\beta$. Call this point $p_3$. Clearly a $k$-space $K$ which contains $p$ will also contain $p_2$ if and only if it contains $p_3$. Hence, also the point $p_2$ indexes a row that is also found in the point-$(k-1)$-space incidence matrix of $\pg(n-1,q)$.
\end{itemize}
Thus we find that,


$$P_n=\begin{pmatrix}
\mathbf{1}^T  & \mathbf{0}^T\\
R_{n-1} & D
\end{pmatrix}\,,$$
where $D$ is an unknown part and $R_{n-1}$ is a matrix of which all rows are also some 
row of of the point-$(k-1)$-space incidence matrix $P_{n-1}$ of $\beta$. 
In fact, each row of $P_{n-1}$ is seen exactly $q$ times in $R_{n-1}$. 
Thus the number of rows of $R_{n-1}$ equals $q\gauss{n}{1}_q=\gauss{n+1}{1}_q-1$, 
while the number of columns equals $\gauss{n}{k}_q$ which is exactly the the number 
of columns of $P_{n-1}$. 

Now suppose that $\chi_\mathcal{L}$ is the characteristic vector of $\mathcal{L}$, 
i.e. $\chi_\mathcal{L} \in \mathrm{Im}(P_n^T)$, then there exists a vector 
$v= (v'\, v_1)^T$, $v' \in \mathbb{R}$ en $v_1 \in \mathbb{R}^{l-1}$, 
$l = \gauss{n+1}{1}_q$ such that

\begin{equation*}
\chi_\mathcal{L} =P_n^T \cdot v =\begin{pmatrix} {\bf 1} & R_{n-1}^T \\ 
{\bf 0} & D^T \end{pmatrix}\cdot \begin{pmatrix} v' \\ v_1 \end{pmatrix} =\begin{pmatrix} v' \cdot {\bf 1} + R_{n-1}^T \cdot v_1 \\ D^T \cdot v_1\end{pmatrix}\,.\\
\end{equation*}
But, due to the definition of $R_{n-1}$, the vector 
$v' \cdot {\bf 1} + R_{n-1}^T \cdot v_1 = P_{n-1}^T \cdot \left(\frac{v'}m\cdot{\bf 1}+v_2\right)$, 
where $m=\gauss{k}{1}_q$ and thus equals the number of points in a $(k-1)$-space. 
Secondly, $v_2$ is a vector of dimension $\gauss{n}{k}_q$, where every 
position $(v_2)_i$ is equal to the sum of $(v_1)_j$, where the rows $(R_{n-1})_{j*}=(P_{n-1})_{i*}$.
Hence the part of $\chi_{\cL}$ representing the $k$-spaces of $\cL$ through $p$, which also represent the
projected $(k-1)$-spaces in $\beta$, belongs to $\mathrm{Im}(P_{n-1}^T)$. So $\cL^p$ is a 
CL $(k-1)$-set in $\beta$.

Assume now that $i>0$, and consider an $(i-1)$-space $\tau$ 
contained in $\pi$. By the induction hypothesis, the set
$\cL^{\tau}$ is a CL $(k-i)$-set in a subspace $\beta = \pg(n-i,q)$,
$\beta$ skew to $\tau$. Now the projection of $\pi$ from $\tau$ 
on $\beta$ is a point $p$. Now the statement follows by applying
again the case $i=0$ on the point $p$.
\end{proof}

Theorem~\ref{DualSubspace} is the dual of Theorem~3.1 in \cite{Me3}. Although these theorems are considered
as folklore by some authors, we included a proof to keep this paper self-contained. A proof of 
Theorem~\ref{DualSubspace} and \cite[Theorem~3.1]{Me3} can also be found as Lemma~6 in \cite{Me5}, formulated
in the language of low degree Boolean functions. Note that Theorem~\ref{DualSubspace} does 
not provide any connection between the parameters of $\cL$ and $\cL^{\pi}$. The following theorem
does. 

\begin{St}\label{th:parameter}
Suppose that $\cL$ is a CL $k$-set with parameter $x$ in $\pg(n,q)$, $n\geq 2k+1$. 
Suppose that $\tau$ is a $(k+1)$-space that contains no $k$-spaces of $\mathcal{L}$.   
If $\pi$ is an $i$-space contained in $\tau$, $i\leq k-2$, then with the 
notation from Theorem \ref{DualSubspace}, $\cL^\pi$ is a CL $(k-i-1)$-set of the same parameter $x$ 
in a subspace $\beta=\pg(n-i-1,q)$.
\end{St}
\begin{proof} 
First assume that $i=0$. Fix an arbitrary point $p\in \tau$ and consider an 
$(n-1)$-space $\beta$, such that $p\not \in \beta$. Then by Theorem~\ref{DualSubspace}, 
$\cL^p$ is a CL $(k-1)$-set in $\beta$. We only have to show that the parameter 
$x^p$ of $\cL^p$ equals $x$. By Theorem~\ref{DrudgeArgum}, we know that
\[
\left| [p]_k\cap \mathcal{L} \right|+ \frac{\gauss{n-1}{k}_q (q^k-1)}{ (q^{k+1}-1)} \left| [\tau]_k\cap \mathcal{L} \right| 
= \gauss{n-1}{k}_q  \left| [p,\tau]_k\cap \mathcal{L} \right| + \frac{q^k-1}{q^n-1} \left| \mathcal{L} \right|\,.
\]
It follows that $\left| [p]_k\cap \mathcal{L} \right|= |\mathcal{L}^p|$ 
and $\left| [\tau]_k\cap \mathcal{L} \right|=\left| [p,\tau]_k\cap \mathcal{L} \right| =0$, 
by choice of $\tau$. Hence,
$$|\mathcal{L}^p|=\frac{q^k-1}{q^n-1}|\mathcal{L}|,$$
or $x^p=x$. \\
Now let $i > 0$. Consider the $i$-dimensional space $\pi \subset \tau$ and choose any point $p \in \pi$. Then 
$\cL^p$ is a CL $(k-1)$-set in $\beta = \pg(n-1,q)$. The subspaces $\pi$ and $\tau$ are projected by $p$ on $\beta$
onto subspaces $\pi'$ of dimension $i-1$, respectively $\tau'$ of dimension $k$ in $\beta$. By induction, 
the set $(\cL^p)^{\pi'}$ is a CL $(k-i-1)$-set of the same parameter $x$ in an $(n-i-1)$-dimensional subspace 
$\beta' \subset \beta$, and clearly, $\cL^{\pi} = (\cL^{p})^{\pi'}$.
\end{proof}

\begin{Gev}\label{th:induction}
Suppose that $\mathcal{L}$ is a CL $k$-set in $\pg(n,q)$, $n\geq 2k+1$, with parameter $x$, then one of the following statements is true:
\begin{enumerate}
\item There exists a CL line class in $\pg(n-k+1,q)$ with parameter $x$.
\item Every $(k+1)$-space $\tau$ contains at least one $k$-space of $\mathcal{L}$.
\end{enumerate}
\end{Gev}
\begin{proof}
    This follows from Theorem \ref{th:parameter}, with $i=k-2$.
\end{proof}

The following lemma shows whenever such skew $(k+1)$-spaces exist.

\begin{Le}\label{Le:RestrictionsOnParameter}
Suppose that $\mathcal{L}$ is a CL $k$-set with parameter $x<\frac{q^{n-k}-1}{q^{k+1}-1}$ in $\pg(n,q)$, then there exists a $(k+1)$-space that contains no element of $\mathcal{L}$.
\end{Le}
\begin{proof}
Suppose that $\mathcal{L}$ is a CL $k$-set with parameter $x<\frac{q^{n-k}-1}{q^{k+1}-1}$. 
Then it is clear that there exists a $k$-space $K\not\in \mathcal{L}$. We want to 
count the number of elements in 
$T:=\{\tau \mid \dim(\tau)=k+1, K\subseteq \tau \text{ containing an element of } \mathcal{L}\}$. 
We do this by double counting the pairs in 
$S:=\{(K_2,\tau)\mid \dim(K_2\cap K)=k-1, \langle K_2, K\rangle= \tau, K_2\in \mathcal{L}, \tau \in T\}$.\\
For a fixed $\tau$, there is at least one $k$-space satisfying the conditions, while 
for every $k$-space $K_2$ there exists only one $(k+1)$-space $\tau$ such that $(K_2,\tau) \in S$. 
Now use that the number of suitable $k$-spaces equals $x\frac{q^{k+1}-1}{q-1}$, 
see Theorem \ref{EquivalenceProj}, to obtain the inequality
$$|T|\cdot 1 \leq x\frac{q^{k+1}-1}{q-1}\cdot 1.$$
In particular, we obtain that $|T|\leq x\frac{q^{k+1}-1}{q-1}$. But if $x<\frac{q^{n-k}-1}{q^{k+1}-1}$, then
$$|T|< \frac{q^{n-k}-1}{q-1},$$
which equals all the $(k+1)$-spaces containing $K$. Thus there exists at 
least one $(k+1)$-space containing $K$ that is not an element of $T$ and 
therefore is skew to the set of $k$-spaces of $\mathcal{L}$.
\end{proof}

\begin{Opm}
Note that if $\mathcal{L}$ consists of all the $k$-spaces in a hyperplane, 
its parameter equals $\frac{q^{n- k}-1}{q^{k+1}-1}$, hence indeed, 
this example is a counterexample to Lemma~\ref{Le:RestrictionsOnParameter} 
if $x$ exceeds its bound. This shows that the bound is sharp.
\end{Opm}

\begin{Gev}\label{cor:nonexis}
The parameter $x<\frac{q^{n-k}-1}{q^{k+1}-1}$ of a CL $k$-set in $\pg(n,q)$, with $n\geq 2k+1$, 
satisfies the same conditions as the parameter of a CL line class in $\pg(n-k+1,q)$.
\end{Gev}
\begin{proof}
Combining Corollary \ref{th:induction} and Lemma \ref{Le:RestrictionsOnParameter}.
\end{proof}

\begin{St}\label{th:modular}
Suppose that $\mathcal{L}$ is a CL $k$-set with parameter $x$ in $\pg(n,q)$, 
such that $n\geq 2k+1$, $n-k+1\geq 7$ and $n-k$ even. Then for every $(k-2)$-dimensional 
subspace $\pi$ that is contained in a $(k+1)$-space that itself contains no elements of $\mathcal{L}$, 
it holds that
$$x(x-1)+2m(m-x)\equiv 0 \mod (q+1),$$
where $m$ denotes the number of $k$-spaces of $\mathcal{L}$ containing $\pi$.\\
\end{St}
\begin{proof}
Suppose that there exists a $(k+1)$-space $\tau$ that contains no 
$k$-spaces of $\mathcal{L}$, and choose an arbitrary $(k-2)$-space 
$\pi$ in $\tau$. By Theorem~\ref{DualSubspace}, the set $\cL^\pi$
is a CL line class in a subspace $\beta = \pg(n-k+1,q)$ skew to $\pi$.
By Theorem~\ref{th:parameter}, $\cL^\pi$ has the same parameter as $\cL$.
By Theorem~\ref{th:projectiveLines} we obtain for any point $p \in \beta$ that
$$x(x-1)+2m(m-x)\equiv 0 \mod (q+1),$$
where $m$ equals the number of lines of $\cL^\pi$ through the point $p$ in $\beta$. 
Each line of $\cL^\pi$ through $p$ in $\beta$ corresponds to a $k$-space through $\pi$ in $\pg(n,q)$.
Hence, we conclude that $m$ equals the number of $k$-spaces of $\cL$ through $\pi$.
\end{proof}
If $x<\frac{q^{n-k}-1}{q^{k+1}-1}$, then it follows from Lemma \ref{Le:RestrictionsOnParameter} 
that there exists a $(k+1)$-space that contains no element of $\mathcal{L}$. This implies 
that the previous result effectively forms a restriction on the parameter $x$.
A second consequence of Corollary \ref{cor:nonexis} is the following.

\begin{Gev}
    Suppose that $\mathcal{L}$ is a CL $k$-set of parameter $x<\frac{q^{n-k}-1}{q^{k+1}-1}$ in $\pg(n,q)$, with $n\geq 2k+1$ and $n\equiv k \mod 2$. Then the parameter $x$ is an integer.
\end{Gev}
\begin{proof}
    This follows from Corollary \ref{cor:nonexis} together with the fact that the parameter of a CL line class in $\pg(n-k+1,q)$ for $n-k+1$ odd is an integer.
\end{proof}

\subsection{The affine case}

From Theorem \ref{CLkSpaceBasic}, it follows that every CL $k$-set in $\ag(n,q)$ is also a CL $k$-set in $\pg(n,q)$ with the same parameter $x$. Hence, Theorem \ref{th:modular} is also valid in the affine case.  To show an even stronger equality, we require the following theorem.

\begin{St}\cite[Theorem 6.16]{Me2}\label{ToLines}
Let $\mathcal{L}$ be a CL $k$-set in AG($n,q$), with $n\geq k+2$. Suppose now that $\mathcal{L}$ has parameter $x$, then $x$ satisfies every condition which holds for Cameron-Liebler line classes in AG($n-k+1,q$).
\end{St}

Combining this with Theorem \ref{th:affine-init} finally results in the following theorem.

\begin{St}\label{ModAf}
Suppose that $\mathcal{L}$ is a CL $k$-set with parameter $x$ in $\ag(n,q)$, with $n-k+1\geq 3 $ and $n-k$ even, then
$x(x-1)\equiv 0 \mod 2(q+1).$
\end{St}

\section*{Acknowledgment}

The authors acknowledge Ferdinand Ihringer for the fruitful discussions on Cameron-Liebler 
$k$-sets and valuable suggestions to improve this manuscript. 


\bibliographystyle{plain}

\end{document}